\documentstyle{article} 
\input amssym.def
\input amssym.tex

\newcommand{\forces}{\Vdash} 

\newcommand{\incomp}{\bot}
\newcommand{\compatible}{\mbox{$\not\!\bot$}}

\newcommand{\V}{{\bf V}} 

\newcommand{\lesdot}{\mathrel{\mathord{<}\!\!\raise 
0.8 pt\hbox{$\scriptstyle\circ$}}}

\newcommand{\cov}{{\rm {\bf cov}}\/}

\newcommand{\con}{{\frak c}}

\newcommand{\cf}{{\rm cf}\/} 
\newcommand{\can}{2^{\textstyle \omega}} 
\newcommand{\fs}{2^{\textstyle <\!\omega}} 
\newcommand{\baire}{\omega^{\textstyle \omega}} 
\newcommand{\iso}{[\omega]^{\textstyle \omega}}

\newcommand{\lh}{{\rm lh}\/} 

\newcommand{\rest}{{\restriction}}

\newcommand{\dom}{{\rm dom}} 
\newcommand{\rng}{{\rm rng}}

\newcommand{\A}{{\cal A}}

\newcommand{\C}{{\cal C}}
\renewcommand{\c}{{\Bbb C}}
\newcommand{\D}{{\cal D}}

\renewcommand{\H}{{\cal H}}  
  
\newcommand{\K}{{\cal K}}

\newcommand{\M}{{\cal M}}

\newcommand{\p}{{\Bbb P}}

\newcommand{\q}{{\Bbb Q}}  

\renewcommand{\r}{{\Bbb R}}

\newcommand{\MA}{{\bf MA}}
\newcommand{\AMA}{{\bf AMA}}
\newcommand{\CH}{{\bf CH}}

\newcommand{\QED}{\hfill\vrule width 6pt height 6pt depth 0pt 
\vspace{0.1in}} 

\newcommand{\Proof}{\noindent{\sc Proof} \hspace{0.2in}} 
\newtheorem{theorem}{Theorem}[section] 
\newtheorem{claim}{Claim}[theorem]

\newtheorem{lemma}[theorem]{Lemma} 
\newtheorem{proposition}[theorem]{Proposition} 
\newtheorem{corollary}[theorem]{Corollary} 
 
\newtheorem{definition}[theorem]{Definition}

\newtheorem{remark}[theorem]{Remark}
\title{Simple forcing notions and Forcing Axioms} 
\author{
{\bf Andrzej Ros{\l}anowski}\thanks{\ \ The research supported by  KBN
(Polish Committee of Scientific Research) grant 1065/P3/93/04}\\
Institute of Mathematics\\
The Hebrew University of Jerusalem\\
91904 Jerusalem, Israel\\
and\\
Mathematical Institute of Wroclaw University\\
50384 Wroclaw, Poland
\and
{\bf Saharon Shelah}\thanks{\ \ The research partially supported by
``Basic Research Foundation'' of the Israel Academy of Sciences and
Humanities. Publication 508}\\
Institute of Mathematics\\
The Hebrew University of Jerusalem\\
91904 Jerusalem, Israel\\
and\\
Department of Mathematics\\
Rutgers University\\
New Brunswick, NJ 08854, USA
}
 
\date\today 
\begin{document}
\maketitle 


\eject

\section{Introduction} 
In the present paper we are interested in simple forcing notions and Forcing
Axioms. A starting point for our investigations was the article \cite{JR1} in
which several problems were posed. We answer some of those problems here.  

In the first section we deal with the problem of adding Cohen reals by simple
forcing notions. Here we interpret {\em simple} as {\em of small size}. We try
to establish as weak as possible versions of Martin Axiom sufficient to
conclude that some forcing notions of size less than the continuum add a
Cohen real. For example we show that $\MA(\sigma\mbox{-centered})$ is enough
to cause that every small $\sigma$-linked forcing notion adds a Cohen real
(see \ref{centimplink}) and $\MA(\mbox{Cohen})$ implies that every small
forcing notion adding an unbounded real adds a Cohen real (see
\ref{unbimpcoh}). A new almost $\baire$-bounding $\sigma$-centered forcing 
notion $\q_{\circledcirc}$ appears naturally here. This forcing notion is
responsible for adding unbounded reals in this sense, that
$\MA(\q_{\circledcirc})$ implies that every small forcing notion adding a new
real adds an unbounded real (see \ref{newimpunb}). 

In the second section we are interested in Anti--Martin Axioms for simple
forcing notions. Here we interpret {\em simple} as {\em nicely definable}. Our
aim is to show the consistency of $\AMA$ for as large as possible class of ccc
forcing notions with large continuum. It has been known that
$\AMA(\mbox{ccc})$ implies $\CH$, but it has been (rightly) expected that
restrictions to regular (simple) forcing notions might help. This is known
under large cardinals assumptions and here we try to eliminate them. We show
that it is consistent that the continuum is large (with no real restrictions)
and $\AMA(\mbox{projective ccc})$ holds true (see \ref{getama}). 

Lastly, in the third section we study the influence of $\MA$ on
$\Sigma^1_3$--absoluteness for some forcing notions. We show that
$\MA_{\omega_1}(\p)$ implies $\Sigma^1_3(\p)$--absoluteness (see
\ref{getabs}).  

\noindent{\bf Notation:} Our notation is rather standard and essentially
compatible with that of \cite{Je} and \cite{BaJu}. However, in forcing
considerations we keep the convention that {\em a stronger condition is the
greater one}.\\ 
For a forcing notion $\p$ and a cardinal $\kappa$ let $\MA_{\kappa}(\p)$ be
the following statement: 
\begin{quotation}
\noindent{\em If $\A_\alpha\subseteq\p$ are maximal antichains in $\p$
(for $\alpha<\kappa$), $p\in\p$

\noindent then there exists a filter $G\subseteq\p$ such that $p\in G$ and
$G\cap \A_\alpha\neq\emptyset$ for all $\alpha<\kappa$.
}
\end{quotation}
For a class $\K$ of forcing notions the sentence $\MA_{\kappa}(\K)$
means $(\forall\p\!\in\!\K)\MA_{\kappa}(\p)$; $\MA_\kappa$ is the sentence
$\MA_\kappa(ccc)$.\\
For a forcing notion $\p$, the canonical $\p$--name for the generic filter on
$\p$ will be called $\Gamma_{\p}$. The incompatibility relation on $\p$ is
denoted by $\incomp_{\p}$ (so $\compatible_{\p}$ means ``compatible'').\\
$\con$ stands for the cardinality of the continuum. For a tree
$T\subseteq\fs$, $[T]$ is the set of all $\omega$--branches through $T$.\\
The family of all sets hereditarily of cardinality $<\chi$ (for a regular
cardinal $\chi$) is denoted by $\H(\chi)$.

\section{Adding a Cohen real}
In this section we obtain several results of the form ``a (weak) version of
$\MA$ implies that small forcing notions (of some type) add Cohen reals''. As
a consequence we answer Problem 5.3 of \cite{JR1} (see \ref{almbound},
\ref{ad5.3con} below). 

\begin{proposition}
\label{addcoh}
Suppose $\p$ is a forcing notion and $\bar{h}$ is a function such that
\begin{enumerate}
\item $\dom(\bar{h})\subseteq \p$, $\rng(\bar{h})\subseteq\fs$,
\item if $p_1,p_2\in\dom(\bar{h})$, $p_1\compatible_{\p} p_2$ then either
$\bar{h}(p_1)\subseteq\bar{h}(p_2)$ or $\bar{h}(p_2)\subseteq\bar{h}(p_1)$, 
\item if $q\in\p$ then there is $\nu_0\in\fs$ such that 
\[(\forall\nu\in\fs, \nu_0\subseteq\nu)(\exists p'\in\dom(\bar{h}))(
p'\compatible_{\p} q\ \&\ \nu\subseteq \bar{h}(p')).\]
\end{enumerate}
Then $\p$ adds a Cohen real.
\end{proposition}

\Proof Though this is immediate, we present the proof fully for reader's
convenience. Let $\bar{h}:\dom(\bar{h})\longrightarrow\fs$ be the function
given by the assumptions. Define a $\p$-name $\dot{c}$ by
\[\forces_{\p}\dot{c}=\bigcup\{\bar{h}(p):
p\in\dom(\bar{h})\cap\Gamma_{\p}\}.\] 
First note that, by the properties of $\bar{h}$, for every filter $G\subseteq
\p$ the set $\{\bar{h}(p): p\in\dom(\bar{h})\cap G\}$ is a chain in
$(\fs,{\subseteq})$. Hence 
\[\forces_{\p}\dot{c}\in 2^{\textstyle{\leq}\omega}.\]
But really $\dot{c}$ is a name for a member of $\can$: suppose not. Then we
have $q\in\p$, $m\in\omega$ such that
\[q\forces_{\p}\dot{c}\in 2^{\textstyle m}.\]
Applying the third property of $\bar{h}$ we get $\nu_0\in\fs$ as there. Let
$\nu\in\fs$, $\nu_0\subseteq\nu$ be such that $\lh(\nu)> m$. We find
$p'\in\dom(\bar{h})$ such that $p'\compatible_{\p} q$ and $\nu\subseteq
\bar{h}(p')$. Thus $p'\forces_{\p}\nu\subseteq\bar{h}(p')\subseteq \dot{c}$,
a contradiction.

To show that 
\[\forces_{\p}\mbox{``}\dot{c}\mbox{ is a Cohen real over {\bf V}''}\]
suppose that we have a closed nowhere dense set $A\subseteq\can$ and a
condition $q\in\p$ such that 
\[q\forces_{\p}\dot{c}\in A.\]
Take $\nu_0\in\fs$ given by condition (3) (for $q$). Since $A$ is nowhere
dense we may choose $\nu\in\fs$ such that $\nu_0\subseteq\nu$ and $[\nu]\cap
A=\emptyset$. By the choice of $\nu_0$, there is a condition $p'\in
\dom(\bar{h})$ such that $p'\compatible_{\p} q$ and $\nu\subseteq
\bar{h}(p')$ (so $p'\forces_{\p} \dot{c}\notin A$), a contradiction. \QED 

\begin{theorem}
\label{centimplink}
Assume $\MA_\kappa(\sigma\mbox{-centered})$. If $\p$ is a $\sigma$-linked
atomless forcing notion of size $\kappa$ then $\p$ adds a Cohen real.
\end{theorem}

\Proof We may assume that the partial order $(\p,\leq)$ is separative, i.e.
\begin{quotation}
\noindent if $p,q\in\p$, $p\nleq_{\p} q$

\noindent then there is $r\in\p$ such that $q\leq_{\p} r$ and $r\incomp_{\p} 
p$. 
\end{quotation}
Of course we may assume that $\p$ is a partial order on a subset of $\can$. We
are going to show that (under our assumptions) there exists a function
$\bar{h}$ as in the assumptions of proposition~\ref{addcoh}. Since $\p$ is 
$\sigma$-linked there are sets $\D_n\subseteq\p$ such that
$\bigcup\limits_{n\in\omega}\D_n =\p$ and any two members of $\D_n$ are
compatible in $\p$ (i.e. each $\D_n$ is linked). Let $N$ be a countable
elementary submodel of $(\H((\beth_7)^+),\in,<^*)$ such that $\p,\langle\D_n:
n\in\omega\rangle,\ldots\in N$.  
\medskip

We define a forcing notion $\r=\r(\p)$: 
\medskip

\noindent{\bf conditions} are pairs $r=\langle h,w\rangle=\langle
h^r,w^r\rangle$ such that 
\begin{description}
\item[(a)] $h$ is a finite function, $\dom(h)\subseteq\p\cap N$,
$\rng(h)\subseteq\fs$,
\item[(b)] if $p_1,p_2\in\dom(h)$ then either $p_1\leq_{\p}p_2$ or
$p_2\leq_{\p}p_1$ or $p_1\incomp_{\p} p_2$,
\item[(c)] if $p_1,p_2\in\dom(h)$, $p_1\leq_{\p} p_2$ then $h(p_1)\subseteq
h(p_2)$,
\item[(d)] $w\in [\p]^{\textstyle {<}\omega}$,
\end{description}

\noindent{\bf the order} is such that $r_1\leq_{\r} r_2$ if and only if
\begin{description}
\item[$(\alpha)$] $h^{r_1}\subseteq h^{r_2}$,
\item[$(\beta)$]  $w^{r_1}\subseteq w^{r_2}$,
\item[$(\gamma)$] {\em if} $q\in w^{r_1}$, $p\in\dom(h^{r_1})$, $p,q$ are
compatible in $\p$ and no $p'\in\dom(h^{r_1})$ satisfies $p\leq_{\p} p'$,
$p\neq p'$, $p'\compatible_{\p} q$

{\em then} the set
\[\begin{array}{ll}
J^{r_1,r_2}_{p,q}\stackrel{\rm def}{=}\{h^{r_2}(p_1): & p\leq p_1\in
\dom(h^{r_2})\ \&\ p_1\compatible_{\p} q\ \&\\
\ & (\forall p_2)(p_1<p_2\in\dom(h^{r_2})\ \Rightarrow\ p_2\incomp_{\p} q)\}
\end{array}
\]
contains a front of $\fs$ above $h^{r_1}(p)$ (i.e. for every $\eta\in\can$
such that $h^{r_1}(p)\subseteq\eta$ there is $k$ with $\eta\restriction k\in
J^{r_1,r_2}_{p,q}$).  
\end{description}

\begin{claim}
\label{cl1}
$(\r,\leq_{\r})$ is a partial order.
\end{claim}

\noindent {\em Proof of the claim:}\hspace{0.15in} The relation $\leq_{\r}$
is reflexive as $J^{r,r}_{p,q}=\{h^r(p)\}$ for all relevant $p,q$. For the
transitivity suppose that $r_1\leq_{\r} r_2$ and $r_2\leq_{\r} r_3$. Clearly
the conditions $(\alpha)$, $(\beta)$ for the pair $r_1,r_3$ are satisfied.
To get condition $(\gamma)$ note that if $\{\nu_0,\ldots,\nu_{k-1}\}$ is a
front in $\fs$ above $\nu$ and $\{\nu^0_0,\ldots,\nu^{l-1}_0\}$ is a front
in $\fs$ above $\nu_0$ then $\{\nu_0^0,\ldots,\nu^{l-1},\nu_1,\ldots,
\nu_{k-1}\}$ is a front above $\nu$. 

\begin{claim}
\label{cl2}
$\r$ is $\sigma$-centered.
\end{claim}

\noindent {\em Proof of the claim:}\hspace{0.15in} Note that if
$r_1,r_2\in\r$, $h=h^{r_1}=h^{r_2}$ then
$\langle h,w^{r_1}\cup w^{r_2}\rangle$ is a common upper bound of $r_1,r_2$.

\begin{claim}
\label{cl3}
Suppose $p\in\p\cap N$, $q\in\p$, $r_0\in\r$ and $m\in\omega$. Then the
following sets are dense in $\r$: 
\begin{enumerate}
\item $I^0_p\stackrel{\rm def}{=}\{r\in\r:\ (\forall q\in
w^r)[p\compatible_{\p} q\ \Rightarrow\ (\exists
p'\in\dom(h^r))(p\leq_{\p}p'\ \&\ p'\compatible_{\p} q)]\}$,
\item $I^1_q\stackrel{\rm def}{=}\{r\in\r:\ q\in w^r\}$,
\item $I^2_{r_0,m}\stackrel{\rm def}{=}\{r\in\r:\ r\incomp_{\r} r_0$\ \
or\ \ for every $q\in w^{r_0}$ and $p\in\dom(h^{r_0})$ such 

\hspace{1.5cm} that $p\compatible_{\p}q\ \&\ (\forall
p'\in\dom(h^{r_0}))([p\leq_{\p} p'\ \&\ p\neq p']\ \Rightarrow\
p'\incomp_{\p} q)$ 

\hspace{1.5cm} and for every $\nu\in 2^{\textstyle m}$ such that
$h^{r_0}(p)\subseteq\nu$ there is 

\hspace{1.5cm} $p''\in\dom(h^r)$ with $p\leq_{\p}p''$, $p''\compatible_{\p} q$
and $\nu\subseteq h^r(p'')\}$.
\end{enumerate}
\end{claim}

\noindent {\em Proof of the claim:}\hspace{0.15in} 1)\ \ \ Assume
$p\in\p\cap N$, $r_0\in\r$. Let $\langle q^l: l<l^*\rangle$ be an
enumeration of $\{q\in w^{r_0}: q\compatible_{\p} p\}$. Choose conditions
$p_l$ (for $l<l^*$) such that 
\begin{enumerate}
\item $p_l\in\p\cap N$
\item for each $p'\in\dom(h^{r_0})$ either $p'\leq_{\p}p_l$ or
$p'\incomp_{\p} p_l$,
\item $p\leq_{\p} p_{l}$,
\item $\langle p_{l}: l<l^*\rangle$ are pairwise incompatible,
\item $p_{l}\compatible_{\p} q^l$.
\end{enumerate}
For this we need the assumption that $\p$ is atomless and $\sigma$--linked.
First take $p_{l}^+\in\p$ such that $p,q^l\leq_{\p} p^+_{l}$ (for $l<l^*$). 
Next we choose $p^{++}_{l}\in\p$, $p^+_{l}\leq_{\p} p^{++}_{l}$ such that
the clauses (2)--(4) are satisfied (remember that $\p$ is atomless and
$\dom(h^{r_0})$ is finite). Let $n_l\in\omega$ be such that $p^{++}_{l}\in
\D_{n_l}$. As $N$ is an elementary submodel of $(\H(\beth_7^+),\in,<^*)$ we
find $\langle p_{l}: l<l^*\rangle\in N$ such that $p_{l}\in\D_{n_l}$ and
the clauses (2)--(4) are satisfied. But now we have (1) too. Moreover this
sequence satisfies (5) since $p_{l}, p^{++}_{l}\in\D_{n_l}$ and the second
condition is stronger than $q^l$ (remember that the sets $\D_{n_l}$ are
linked). 

Define $h^r$ by
\[\dom(h^r)=\dom(h^{r_0})\cup\{p_l: l<l^*\},\ \ h^{r_0}\subseteq h^r\ \ \mbox{
and}\] 
\[h^r(p_l)=\bigcup\{h^{r_0}(p'): p'\leq_{\p}p_l\ \&\ p'\in\dom(h^{r_0})\}.\]
First note that all conditions $p'\in\dom(h^{r_0})$ satisfying $p'\leq_{\p}
p_l$ are compatible in $\p$ and hence (by {\bf (b)} for $r_0$) they are
pairwise comparable and thus (by {\bf (c)} for $r_0$) the set
\[\{h^{r_0}(p'): p'\leq_{\p} p\ \&\ p'\in\dom(h^{r_0})\}\]
is a (finite) chain in $(\fs,{\subseteq})$.  Hence $h^r(p_l)\in\fs$ (actually
$h^r(p_l)=h^{r_0}(p^*_l)$ for the $\leq_{\p}$-maximal
$p^*_l\in\dom(h^{r_0})$ such that $p^*_l\leq_{\p} p_l$; if there is no such
$p^*_l$ then $h^r(p_l)=\langle\rangle$). Consequently $h^r$ satisfies {\bf
(a)}. One can easily check that $h^r$ satisfies conditions {\bf (b)}, {\bf (c)}
too and thus $r=\langle h^r,w^{r_0}\rangle\in I^0_p$.  Conditions
$(\alpha),(\beta)$ for the pair $r_0,r$ are clear. To check the clause
$(\gamma)$ suppose that $p'\in\dom(h^{r_0})$, $q\in w^{r_0}$ are relevant
for it. If for each $l<l^*$ either $p_l\incomp_{\p} q$ or $p'\nleq_{\p} p_l$
(which implies $p'\incomp_{\p} p_l$) then $J^{r_0,r}_{p',q}=\{h^{r_0}(p')\}$
as the property of $p'$ there is preserved. Otherwise
$J^{r_0,r}_{p',q}=\{h^{r}(p_l): l<l^*, p_l\compatible_{\p} q, p'\leq_{\p}
p_l\}$. But due to condition {\bf (b)} for $r_0$ we 
have that each condition from $\dom(h^{r_0})$ weaker than any $p_l$ such
that $p_l\compatible_{\p} q$, $p'\leq_{\p} p_l$ is weaker than
$p'$. Consequently $h^r(p_l)=h^{r_0}(p')$ for all relevant $p_l$ and we get
$r_0\leq_{\r} r$. 
\medskip

\noindent 2)\ \ \ Let $q\in\p$, $r\in\r$. Take $\langle
h^r,w^r\cup\{q\}\rangle$; easily it is a condition in $I_q$ stronger than
$r$.  
\medskip

\noindent 3)\ \ \ Assume $r_0\in\r$, $m\in\omega$. Let $r\in\r$. If $r_0,r$
are incompatible in $\r$ then $r\in I^2_{r_0,m}$ and we are done. So we may
assume that $r_0\leq_{\r} r$. 

Let $\langle (q^l,p^l,\nu^l):l<l^*\rangle$ list all triples $(q,p,\nu)$ such
that 
\begin{quotation}
\noindent $q\in w^r$, $p\in\dom(h^r)$, $q\compatible_{\p}p$,
$h^r(p)\subseteq\nu\in 2^{\textstyle m}$ and there is no $p'\in\dom(h^r)$
with $p<_{\p} p'$, $q\compatible_{\p} p'$. 
\end{quotation}
(It is possible that $l^*=0$, e.g. if $m$ is too small.)

\noindent Now choose conditions $p^*_l$ such that
\begin{enumerate}
\item $p^*_l\in\p\cap N$,
\item for each $p\in\dom(h^{r})$ either $p\leq_{\p}p^*_l$ or $p\incomp_{\p}
p^*_l$, 
\item $p^l\leq_{\p} p_{l}^*$,
\item $\langle p_{l}^*: l<l^*\rangle$ are pairwise incompatible,
\item $p_{l}^*\compatible_{\p} q^l$.
\end{enumerate}
For this we follow exactly the lines of the respective part of the proof of
1) (so this is another place we use the assumptions on $\p$). 

Next define $h^{r_1}=h^r\cup\{(p^*_l,\nu^l): l< l^*\}$, $w^{r_1}=w^r$,
$r_1=\langle h^{r_1}, w^{r_1}\rangle$. Similarly as in 1) one checks that
$r_1\in \r$. 

The condition $r_1$ is stronger than $r$: clauses $(\alpha)$, $(\beta)$ are
clear. For $(\gamma)$ suppose that $q\in w^r$, $p\in\dom(h^r)$ are relevant
for this clause. If $m\geq\lh(h^r(p))$ then each $\nu\in 2^{\textstyle m}$
extending $h^r(p)$ appears as $\nu^l=h^{r_1}(p^*_l)$ for some $l<l^*$ such
that $q^l=q$, $p^l=p$. Hence $J^{r,r_1}_{p,q}$ contains a front above
$h^r(p)$. If $m<\lh(h^r(p))$ then the pair $(p,q)$ does not appear as
$(p^l,q^l)$. Note that for each $l<l^*$, if $q^l\compatible_{\p} p$ then
$p^l\incomp_{\p} p$ (as $p^l$ cannot be stronger than $p$ since
$\lh(h^{r_1}(p^l))\leq m$) and hence $p^*_l\incomp_{\p} p$. If
$q^l\incomp_{\p} p$ then we get the same conclusion (though $p^l$ might be
weaker than $p$, the demands (5), (2) of the choice of $p^*_l$ imply that
$p^*_l\incomp_{\p} p$). Consequently the ``maximality'' property of $p$ is
preserved in $\dom(h^{r_1})$ and $J^{r,r_1}_{p,q}=\{h^r(p)\}$.

To prove that $r_1\in I^2_{r_0,m}$ suppose that $q\in w^{r_0}$ and
$p\in\dom(h^{r_0})$ is maximal (in $\dom(h^{r_0})$) compatible with $q$. Let
$\nu\in 2^{\textstyle m}$ extend $h^{r_0}(p)$. Since $r_0\leq_{\r} r$ we
find $p'\in\dom(h^r)$ stronger than $p$, maximal (in $\dom(h^r)$) compatible
with $q$ and such that $\nu$, $h^r(p')$ are comparable (by condition
$(\gamma)$). If $\nu\subseteq h^r(p')$ then we are done. So suppose
$h^r(p')\subseteq\nu$. Then for some $l<l^*$ we have $q=q^l$, $p'=p^l$ and
$\nu=\nu^l$. By the choice of $p^l$ and the definition of $h^{r_1}(p_l)$ we
get 
\[p\leq_{\p}p'\leq_{\p}p^*_l\quad \&\quad p^*_l\compatible_{\p} q_l=q\quad
\&\quad \nu=\nu_l=h^{r_1}(p^*_l).\]  
The claim is proved.
\bigskip

Since we have assumed $\MA_\kappa(\sigma\mbox{-centered})$ we find a filter
$H\subseteq\r$ such that  
\begin{description}
\item[$(\oplus_0)$] $H\cap \{r\in\r: r\incomp_{\r}r_0\mbox{ or }(r_0\geq r\
\&\ r\in I^0_p)\}\neq\emptyset$ for $p\in\p\cap N$, $r_0\in\r$,
\item[$(\oplus_1)$] $H\cap I^1_q\neq\emptyset$ for $q\in\p$ and
\item[$(\oplus_2)$] $H\cap I^2_{r_0,m}\neq\emptyset$ for $r_0\in\r$,
$m\in\omega$.
\end{description}
Put $\bar{h}=\bigcup\{h^r:\ r\in H\}$. Clearly $\bar{h}$ is a function from
a subset of $\p\cap N$ to $\fs$. Conditions {\bf (b)}, {\bf (c)} imply that
$\bar{h}$ satisfies the second requirement of the assumptions of
\ref{addcoh}. 

Suppose now that $q\in\p$. Take $p\in\p\cap N$ compatible with $q$ and choose
$r_0\in H\cap I^0_p\cap I^1_q$ (so $q\in w^{r_0}$). Next take
$p^*\in\dom(h^{r_0})$ such that  
\[p\leq_{\p} p^*\ \&\ p^*\compatible_{\p} q\ \&\ (\forall p'\in\dom(h_0))
([p^*\leq_{\p} p'\ \&\ p^*\neq p']\ \Rightarrow\ p'\incomp_{\p} q).\]
Assume that $h_{r_0}(p^*)\subseteq\nu\in 2^{\textstyle m}$. 

\noindent By $(\oplus_2)$ we find $r\in H\cap I^2_{r_0,m}$. As $r_0,r\in
H$, $H$ is a filter, we cannot have $r\incomp_{\r} r_0$. Consequently ``the
second part'' of the definition of $I^2_{r_0,m}$ applies to $r$. Looking
at this definition (with $p^*$ as $p$ there) we see that there is
$p'\in\dom(h^r)$ with 
\[p^*\leq_{\p} p'\ \&\ p'\compatible_{\p} q\ \&\ \nu\subseteq h^r(p').\]
So $\nu_0=\bar{h}(p^*)$ is as required in 3). Applying \ref{addcoh} we finish
the proof of the theorem. \QED
\bigskip

\begin{remark}
{\em
Of course, what we have shown in \ref{centimplink} is that
$\MA_\kappa(\r(\p))$ implies $\p$ adds a Cohen real, provided $\p$ is
atomless $\sigma$--linked of size $\kappa$. 
}
\end{remark}

\begin{corollary}
Assume $\MA_\kappa$. If $\p$ is a ccc atomless forcing notion of size
$\kappa$ then $\p$ adds a Cohen real.\QED
\end{corollary}

\begin{proposition}
\label{unbounded}
Let $\p$ be a ccc forcing notion. Then the following conditions are
equivalent:
\begin{description}
\item[(a)] $\forces_{\p}$``there is an unbounded real in $\baire$ over
$\V$''
\item[(b)] there exists a sequence $\langle \A_n: n\in\omega\rangle$ of
maximal antichains of $\p$ such that
\begin{enumerate}
\item $\A_0=\{0\}$,
\item $(\forall n\in\omega)(\forall p\in\A_{n+1})(\exists q\in\A_n)(
q\leq_{\p} p)$,
\item $(\forall n\in\omega)(\forall p\in\A_n)(\|\{q\in\A_{n+1}:
p\leq_{\p}q\}\|=\omega)$,
\item $(\forall q\in\p)(\exists n\in\omega)(\|\{p\in\A_n: p\compatible_{\p}
q\}\|=\omega)$. 
\end{enumerate}
\end{description}
\end{proposition}

\Proof Easy, left for the reader. \QED

\begin{theorem}
\label{unbimpcoh}
Suppose that $\p$ is a ccc forcing notion such that $\|\p\|<\cov(\M)$ (i.e.
unions of $\|\p\|$ many meager sets are meager) and
\[\forces_{\p}\mbox{``there is an unbounded real over $\V$''}.\]
Then
\[\forces_{\p}\mbox{``there is a Cohen real over $\V$''}.\]
\end{theorem}

\Proof We are going to apply proposition~\ref{addcoh} and for this we will
construct a function $\bar{h}$ satisfying (1)--(3) of \ref{addcoh}.\\
Let $\langle \A_n: n\in\omega\rangle$ be a sequence of maximal antichains of
$\p$ given by {\bf (b)} of proposition~\ref{unbounded}. Take a countable
elementary submodel $N$ of $(\H(\beth_7^+),{\in},{<^*})$ such that $\p,\langle
\A_n:n\in\omega\rangle,\ldots\in N$. Consider the following partial order:
\medskip

\noindent {\bf conditions} are finite functions $h$ such that 
\begin{description}
\item[a.] $\dom(h)\subseteq\bigcup\limits_{n\in\omega}\A_n$,
$\rng(h)\subseteq\fs$, 
\item[b.] if $p_1,p_2\in\dom(h)$, $p_1\leq_{\p} p_2$ then $h(p_1)\subseteq
h(p_2)$,
\end{description}

\noindent{\bf the order} is the inclusion; $h_1\leq_{\c} h_2$ iff
$h_1\subseteq h_2$.  
\medskip

\noindent Clearly $\c$ is (isomorphic to) the Cohen forcing notion.
\begin{claim}
\label{cl4}
Let $p\in\bigcup\limits_{n\in\omega}\A_n$, $h_0\in\c$, $q\in\p$,
$m\in\omega$. Then the following sets are dense in $\c$:
\begin{enumerate}
\item $J^0_p\stackrel{\rm def}{=}\{h\in\c: p\in\dom(h)\}$,
\item $J^1_{q,m,h_0}\stackrel{\rm def}{=}\{h\in\c: h\incomp_{\c} h_0$ or
for every $p\in\dom(h_0)$ such that for some

\hspace{1.5cm} $n\in\omega$, $p\in\A_n$ and the set $\{p'\in\A_{n+1}:
p\leq_{\p} p'\ \&\ p'\compatible_{\p} q\}$ is 

\hspace{1.5cm} infinite we have:\ \ \ \ \ for every $\nu\in 2^{\textstyle
m}$ extending $h_0(p)$ 

\hspace{1.5cm} there is $p'\in\dom(h)$ with $p\leq_{\p} p'$,
$p'\compatible_{\p} q$ and $h(p')=\nu\ \}$
\end{enumerate}
\end{claim}

\noindent {\em Proof of the claim:}\hspace{0.15in} 1)\ \ \ Assume
$p\in\bigcup\limits_{n\in\omega}\A_n$, $h\in\c$. Extend $h$ to $h'$ by
putting 
\[h'(p)=\bigcup\{h(p'): p'\in\dom(h)\ \&\ p'\leq_{\p} p\}.\]
Easily this $h'$ satisfies $h'\in J^0_p$, $h\leq_{\c} h'$.
\medskip 

\noindent 2)\ \ \ Suppose that $q\in\p$, $m\in\omega$, $h_0\in\c$,
$h\in\c$. We may assume that $h_0\leq_{\c} h$. Let $\langle
(p^l,n^l,\nu^l): l<l^*\rangle$ enumerate all triples $p\in\dom(h_0)$,
$n\in\omega$, $\nu\in 2^{\textstyle m}$ such that 
\begin{description}
\item[$(\alpha)$] $p\in\A_n$ and the set $\{p'\in\A_{n+1}: p\leq_{\p} p'\ \&\
p'\compatible_{\p} q\}$ is infinite
\item[$(\beta)$]  $h_0(p)\subseteq\nu$.
\end{description}
Next (using $(\alpha)$ above) choose $p^*_l\in\A_{n^l+1}$ such that
\begin{enumerate}
\item $p^l\leq_{\p} p^*_l$
\item $\langle p^*_l: l<l^*\rangle$ are pairwise incompatible
\item for each $p\in\dom(h)$, $l<l^*$ either $p\leq_{\p} p^*_l$ or
$p\incomp_{\p} p^*_l$.
\end{enumerate}
Now put $\dom(h')=\dom(h)\cup\{p^*_l: l<l^*\}$, $h'(p^*_l)=\nu^l$ and
$h'\restriction\dom(h)=h$. Easily $h'\in\c$, $h\leq_{\c} h'$ and $h'\in
J^1_{q,m,h_0}$. This finishes the claim.
\bigskip

Since $\|\p\|<\cov(\M)$ we find a filter $H\subseteq\c$ such that $H\cap
J^0_p\neq\emptyset$ and $H\cap J^1_{q,m,h_0}\neq\emptyset$ for all $q\in\p$,
$m\in\omega$, $h_0\in\c$ and $p\in\bigcup\limits_{n\in\omega}\A_n$. Put
$\bar{h}=\bigcup H$. Then clearly
$\bar{h}:\bigcup\limits_{n\in\omega}\A_n\longrightarrow\fs$ is a function 
satisfying the requirements (1), (2) of \ref{addcoh}. To check the third
condition there suppose $q\in\p$. Take $n\in\omega$ and $p^*\in\A_n$ such
that the set $\{p'\in\A_{n+1}: p^*\leq_{\p} p'\ \&\ p'\compatible_{\p} q\}$
is infinite (possible by the choice of the $\A_k$'s). Since $H\cap
J^0_{p^*}\neq\emptyset$ we find a condition $h_0\in H$ such that
$p^*\in\dom(h_0)$. Suppose that $\nu\in\fs$, $\bar{h}(p^*)\subseteq\nu$ and
let $m=\lh(\nu)$. Take $h_1\in H\cap J^1_{q,m,h_0}$. Since $h_0,h_1$ cannot be
incompatible, $p^*\in\dom(h_0)$, $h(p^*)\subseteq\nu\in 2^{\textstyle m}$ we
find $p'\in\dom(h_1)$ such that $p^*\leq p'$, $p'\compatible_{\p} q$ and
$h_1(p')=\nu$. Since $h_1(p')=\bar{h}(p')$ we conclude that
$\nu_0=\bar{h}(p^*)$ is as required in (3) of \ref{addcoh} for $q$. The
theorem is proved. \QED

\begin{definition}
A forcing notion $\p$ is {\em almost $\baire$-bounding} if
\begin{quotation}
\noindent for each $\p$-name $\dot{f}$ for an element of $\baire$ and a
condition $p\in\p$ there is $g\in\baire\cap\V$ such that for every
$X\in\iso\cap\V$:
\[(\exists p'\geq_{\p} p)(p'\forces_{\p}(\exists^\infty n\in X)(\dot{f}(n)<
g(n))).\] 
\end{quotation}
\end{definition}

\begin{lemma}
\label{product}
\begin{enumerate}
\item Suppose that $\p$ is a ccc forcing notion such that for every integer $n$
the product forcing notion $\p^n$ does not add unbounded real and satisfies
the ccc. Then the $\omega$-product $\p^\omega$ with finite support is almost
$\baire$-bounding and satisfies the ccc. 
\item Finite support iteration of ccc almost $\baire$-bounding forcing
notions does not add a dominating real. 
\end{enumerate}
\end{lemma}

\Proof 1)\ \ \ Suppose that for each $n\in\omega$ the product forcing notion
$\p^n$ satisfies the ccc and does not add unbounded reals. By \cite[23.11]{Je}
we know that then $\p^\omega$ satisfies the ccc. We have to show that
$\p^\omega$ is almost $\baire$-bounding. Let $\dot{f}$ be a $\p^\omega$--name
for a function in $\baire$. For each $n,k\in\omega$ choose a maximal
antichain $\A^n_k$ of $\p^n$ and mappings $\varphi^n_k:\A^n_k\longrightarrow
\p^\omega$ and $g^n_k:\A^n_k\longrightarrow\omega$ such that
\[(\forall q\in \A^n_k)(\varphi^n_k(q)\rest n=q\ \ \&\ \ \varphi^n_k(q)
\forces_{\p^\omega}\dot{f}(k)=g^n_k(q))\]
(possible as $\p^n\lesdot\p^\omega$). Thus, for each $n\in\omega$, we have a
$\p^n$--name $\dot{g}^n$ for a function in $\baire$ defined by
\[(\forall k\in\omega)(\forall q\in\A^n_k)(q\forces_{\p^n}\dot{g}^n(k)=g^n_k(
q)).\]
Since $\p^n$ does not add unbounded reals and satisfies the ccc we find a
function $g_n\in\baire$ such that
\[\forces_{\p^n}(\exists m\in\omega)(\forall k\geq m)(\dot{g}^n(k)<g_n(k)).\]
Take $g\in\baire$ such that $(\forall n\in\omega)(\exists m\in\omega)(\forall
k\geq m)(g_n(k)<g(k))$. We claim that 
\[\forces_{\p^\omega}(\forall X\in\iso\cap\V)(\exists^\infty k\in
X)(\dot{f}(k)<g(k)).\]
To this end suppose that $X\in\iso$, $p\in\p^\omega$ and $N\in\omega$. Take
$n$ such that $p\in\p^n$ and look at the function $g_n$. By its choice we find
a condition $p'\in\p^n$ stronger than $p$ and an integer $m_0$ such that $p'
\forces_{\p^n}(\forall k\geq m_0)(\dot{g}^n(k)<g_n(k))$. By the choice of $g$
we find $m_1\in\omega$ such that $(\forall k\geq m_1)(g_n(k)<g(k))$. Let $k\in
X$ be such that $k>m_0+m_1+N$. Since $\A^n_k$ is a maximal antichain of
$\p^n$ we may take a condition $q\in\A^n_k$ compatible with $p'$. Let $p''$ be
a common upper bound of $p'$ and $\varphi^n_k(q)$ in $\p^\omega$. Then ($p''$
is stronger than $p$ and) 
\[p''\forces_{\p^\omega}\dot{f}(k)=g^n_k(q)=\dot{g}^n(k)<g_n(k)<g(k)\]
(remember $k$ is above $m_0,m_1$). Since $k\in X$ is greater than $N$ we
finish by standard density arguments.

\noindent 2)\ \ \ See \cite[Ch VI, 3.6+3.17]{Sh:f} or \cite[6.5.3]{BaJu}. \QED

\begin{theorem}
\label{almbound}
Assume $\MA_\kappa(\mbox{ccc}\ \&\ \mbox{almost $\baire$-bounding})$. Then
every atomless ccc forcing notion of size $\leq\kappa$ adds a Cohen real.
\end{theorem}

\Proof We assume of course that $\kappa\geq\aleph_1$. Let $\p$ be a ccc
forcing notion, $\|\p\|\leq\kappa$. If $\p$ adds an unbounded real then
theorem~\ref{unbimpcoh} applies (note that the Cohen forcing notion is almost
$\baire$-bounding, so our assumption implies $\kappa<\cov(\M)$). Thus to
finish the proof we need to show that $\p$ adds an unbounded real. This fact
is done by the two claims below.

\begin{claim}
\label{cl5} Assume $\MA_\kappa(\mbox{ccc}\ \&\ \baire\mbox{-bounding})$.
Suppose that $\p$  is a ccc forcing notion which adds no unbounded real (i.e.
it is $\baire$-bounding). Then for every $n\in\omega$ the product forcing
notion $\p^n$ adds no unbounded real and satisfies the ccc. 
\end{claim}

\noindent{\em Proof of the claim:}\hspace{0.15in} As $\MA_\kappa(\mbox{ccc}\
\&\ \baire\mbox{-bounding})$ applies to $\p$, this forcing notion
has the Knaster property (strong ccc) and consequently all powers of it
satisfy the ccc. What might fail is not adding unbounded reals. So suppose
that $n$ is the first such that 
\[\forces_{\p^n}\mbox{``there is an unbounded real over $\V$''}.\]
Clearly $n>1$. By proposition~\ref{unbounded} we find maximal antichains
$\A_k\subseteq\p^n$ (for $k<\omega$) satisfying conditions (1)--(4) of
clause {\bf (b)} there.

We may think that $\p$ is an ordering on $\kappa$. Let $N$ be an elementary
submodel of $(\H(\beth_7^+),{\in},{<^*})$ such that
\[\p,\leq_{\p},\langle\A_k:k\in\omega\rangle,\ldots\in N,\quad
\kappa+1\subseteq N\quad\mbox{ and }\quad \|N\|=\kappa.\]
Let $\pi:N\longrightarrow M$ be the Mostowski collapse of $N$, $M$ a
transitive set. Note that $\pi(\p)=\p$, $\pi(\A_k)=\A_k$ etc. Since
$\p^{n-1}$ is ccc and adds no unbounded real we may apply our restricted
version of $\MA_\kappa$ to it and get an $M$-generic filter
$H\subseteq\p^{n-1}$ in $\V$. (Note that if $\A\subseteq\p^{n-1}$, $\A\in M$
then $M\models$``$\A$ is a maximal antichain of $\p^{n-1}$''\ \ iff\ \ $\A$ is
really a maximal antichain of $\p^{n-1}$.) Let 
\[\A^H_k\stackrel{\rm def}{=}\{p\in\p: (\exists \bar{p}\in H)((\bar{p},p)\in
\A_k)\}\in M[H].\] 
Then 
\[M[H]\models\mbox{``$\A^H_k$ is a maximal antichain of $\p$''}\]
and easily the same holds in $\V$. As $\p$ adds no unbounded real, by
\ref{unbounded} we find $p\in\p$ such that 
\[(\forall k\in\omega)(\|\{p'\in\A^H_k: p\compatible_{\p} p'\}\|<\omega)\]
and thus
\[(\forall k\in\omega)(\|\{(\bar{p}',p')\in\A_k: \bar{p}'\in H\ \&\
p\compatible_{\p} p'\}\|<\omega).\]
Since $\p^{n-1}$ adds no unbounded real (and this is true in $M$ too) we
find finite sets $A_k\subseteq\A_k$ (for $k\in\omega$) and a condition
$\bar{p}\in\p^{n-1}$ such that for each $k\in\omega$
\[M\models\ \bar{p}\forces_{\p^{n-1}}\{(\bar{p}',p')\in\A_k: \bar{p}'\in
\Gamma_{\p^{n-1}} \ \&\ p\compatible_{\p} p'\}\subseteq A_k.\]
This means that if $(\bar{p}',p')\in\A_k\setminus A_k$ then either
$\bar{p}\incomp_{\p^{n-1}}\bar{p}'$ or $p\incomp_{\p} p'$. Hence the
condition $(\bar{p},p)\in\p^n$ is a counterexample to the fourth property of
$\langle\A_k: k\in\omega\rangle$. The claim is proved.
\bigskip

It follows from \ref{cl5} and \ref{product} that (under our assumptions) $\p$
is $\sigma$-centered. So now we may use the following claim.
\begin{claim}
\label{cl8}
Every $\sigma$-centered atomless forcing notion adds an unbounded real.
\end{claim}

\noindent{\em Proof of the claim:}\hspace{0.15in} Folklore; see e.g. 5.2 of
\cite{JR1}. \QED

\begin{corollary}
\label{ad5.3con}
It is consistent that $\con>\aleph_1$, every atomless ccc forcing notion of
the size $<\con$ adds a Cohen real but $\MA_{\omega_1}(\mbox{ccc})$ fails.\QED 
\end{corollary}

As we saw in \ref{unbimpcoh}, if we assume a small part of $\MA_\kappa$ then
each forcing notion adding an unbounded real adds a Cohen real, provided the
size of the forcing is at most $\kappa$. Therefore it is natural to look for
requirements implying that small forcing notions add unbounded reals. The
main part of the proof of \ref{almbound} was to show that
$\MA_\kappa(\mbox{ccc}\ \&\ \mbox{almost $\baire$-bounding})$ is such a
condition. It occurs however, that we need much less for this. As in
\ref{unbimpcoh} the crucial role was played by the Cohen forcing, here we
naturally arrive to the forcing notion defined below.

\begin{definition}
We define a forcing notion $\q_\circledcirc$:
\medskip

\noindent{\bf conditions} are pairs $\langle a,w\rangle$ such that $w\in
[\can]^{\textstyle {<}\omega}$ and $a\in [\fs]^{\textstyle {<}\omega}$,
\medskip

\noindent{\bf the order} is defined by:\ \ \ $\langle
a_0,w_0\rangle\leq_{\q_\circledcirc}\langle a_1,w_1\rangle$ if and only if

\noindent $a_0\subseteq a_1$, $w_0\subseteq w_1$ and $(\forall\eta\in
w_0)(\forall l\in\omega)(\eta\restriction l\in a_1\ \Rightarrow\
\eta\restriction l\in a_0)$.
\end{definition}

\begin{lemma}
\label{forcingproperties}
\begin{enumerate}
\item $\q_\circledcirc$ is an almost $\baire$-bounding $\sigma$-centered
partial order.  
\item Let $\dot{A}$ be the $\q_\circledcirc$--name for a subset of $\fs$ given
by 
\[\forces_{\q_\circledcirc}\dot{A}=\bigcup\{a: (\exists w)(\langle
a,w\rangle\in \Gamma_{\q_\circledcirc})\}.\]
Then 
\begin{description}
\item[$(\alpha)$] $\forces_{\q_\circledcirc}(\forall\eta\in\can\cap\V)
(\forall^\infty n\in\omega)(\eta\rest n\notin\dot{A})$
\item[$(\beta)$] $\forces_{\q_\circledcirc}$``if $T\subseteq\fs$ is a perfect
tree from the ground model

\hspace{0.7in}then $(\exists^\infty n\in\omega)(T\cap 2^{\textstyle n}\cap
\dot{A}\neq \emptyset)$''. 
\end{description}
\end{enumerate}
\end{lemma}

\Proof 1)\ \ \ Clearly if $a_0=a_1$,
$\langle a_0,w_0\rangle,\langle a_1,w_1\rangle\in\q_\circledcirc$ then
$\langle  a_0,w_0\cup w_1\rangle\in\q_\circledcirc$ is a common upper bound of
$\langle a_0,w_0\rangle$, $\langle a_1,w_1\rangle$. This implies that
$\q_\circledcirc$ is $\sigma$-centered. Next note that 
\begin{quotation}
\noindent $\langle a_0,w_0\rangle\incomp_{\q_\circledcirc} \langle
a_1,w_1\rangle$ if and only if  

\noindent either there are $\eta\in w_0$, $l\in\omega$ such that
$\eta\restriction l\in a_1\setminus a_0$

\noindent or the symmetrical condition holds (interchanging 0 and 1).
\end{quotation}
Hence if $a_0,a_1\subseteq 2^{\textstyle \leq l_0}$, $\{\eta\restriction
l_0:\eta\in w_1\}=\{\eta\restriction l_0: \eta\in w_2\}$ then 
\[\langle a_0,w_0\rangle\incomp_{\q_\circledcirc}\langle a_1,w_1\rangle\
\mbox{ iff }\ 
\langle a_0,w_0\rangle\incomp_{\q_\circledcirc}\langle a_1,w_2\rangle.\] 
Since the product space $(\can)^n$ is compact we may conclude that
\begin{quotation}
\noindent {\em if $\A\subseteq\q_\circledcirc$ is a maximal antichain,
$n\in\omega$, $a\in [\fs]^{\textstyle{<}\omega}$ 

\noindent then there is a finite set $A=A^{a,n}\subseteq\A$ such that  for
every $w\subseteq\can$, $\|w\|=n$ there is $r\in A$ with $\langle
a,w\rangle\compatible_{\q_\circledcirc} r$. 
}
\end{quotation}
The above property easily implies that $\q_\circledcirc$ is almost
$\baire$-bounding: suppose that $\dot{h}$ is a $\q_\circledcirc$-name for an
element of $\baire$. For each $k\in\omega$ fix a maximal antichain $\A_k$ such
that each member of $\A_k$ decides the value of $\dot{h}(k)$. For
$k,n\in\omega$ and $a\in [\fs]^{{<}\omega}$ choose a finite set
$A^{a,n,k}\subseteq\A_k$ with the property stated above. Finally put 
\[g(k)=1+\max\{l\in\omega: (\exists a\subseteq 2^{\textstyle{\leq}k})
(\exists n\leq k)(\exists r\in A^{a,n,k})(r\forces_{\q_\circledcirc}
\dot{h}(k)= l)\}.\]   
To show that the function $g$ works for $\dot{h}$ (for the definition of
almost $\baire$-bounding) suppose that $X\in\iso$. Assume that 
\[r_0\forces_{\q_\circledcirc}(\forall^\infty n\in X)(g(n)\leq\dot{h}(n)),\]
so we have $r_1$ and $k$ such that
\[r_1\forces_{\q_\circledcirc}(\forall n>k)(n\in X\ \Rightarrow\ g(n)\leq
\dot{h}(n)).\]
Now take $k^*\in X$ such that $k^*>k$ and if $r_1=\langle a,w\rangle$ then
$a\subseteq 2^{\textstyle{\leq}k^*}$, $\|w\|=n\leq
k^*$. By the definition of $A^{a,n,k^*}$ we find $r\in A^{a,n,k^*}$
compatible with $r_1$. But each member of $A^{a,n,k^*}$ forces that
$\dot{h}(k^*)<g(k^*)$, a contradiction. 

\noindent 2) Straightforward. \QED

\begin{theorem}
\label{newimpunb}
Assume $\MA_\kappa(\q_\circledcirc)$.
Suppose that $\p$ is a forcing notion such that $\|\p\|\leq\kappa$ and
$\forces_{\p}\can\cap\V\neq\can$ (i.e. the corresponding complete
Boolean algebra ${\rm RO}(\p)$ is not $(\omega,\omega)$-distributive). Then
$\p$ adds an unbounded real.
\end{theorem}

\Proof Since $\p$ adds new reals we can find a $\p$--name $\dot{r}$ for an
element of $\can$ such that $\forces_{\p}\dot{r}\notin\V$. For a condition
$q\in\p$ let
\[T^q\stackrel{\rm def}{=}\{\nu\in\fs: q\not\forces_{\p}\nu\nsubseteq\dot{r}
\}.\] 
By our assumptions on $\dot{r}$ we know that each $T^q$ is a perfect tree in
$\fs$. Next fix $\eta_q\in [T^q]$ (for $q\in \p$). Since we have assumed
$\MA_\kappa(\q_\circledcirc)$ we may apply lemma~\ref{forcingproperties} to
find a set $A\subseteq\fs$ such that for each $q\in\p$:
\begin{description}
\item[$(\alpha)$] $(\forall^\infty n\in\omega)(\eta_q\rest n\notin A)$ and
\item[$(\beta)$]  $(\exists^\infty n\in\omega)(T^q\cap 2^{\textstyle n}\cap
A\neq\emptyset)$. 
\end{description}
Now define a $\p$-name $\dot{K}$ for a subset of $\omega$ by:
\[\forces_{\p}\dot{K}=\{n\in\omega: \dot{r}\rest n\in A\}.\]
First note that $\dot{K}$ is a $\p$-name for an infinite subset of $\omega$:
Why? Suppose that $q\in\p$ and $N\in\omega$. By the property $(\beta)$ of $A$
we find $\nu\in A\cap T^q$ such that $\lh(\nu)>N$. Then we have a condition
$p_\nu\geq q$ which forces ``$\nu\subseteq\dot{r}$'' and thus
$p_\nu\forces_{\p}\lh(\nu)\in \dot{K}$.

\noindent Suppose now that $q\in\p$, $g\in\baire$ is an increasing function
and $N_0\in\omega$. Take $N_1>N_0$ such that $(\forall n\geq N_1)(\eta_q\rest 
n\notin A)$ and a condition $p_{\eta_q\rest g(N_1)}$ such that $q\leq_{\p}
p_{\eta_q\rest g(N_1)}$ and $p_{\eta_q\rest g(N_1)}\forces_{\p} \eta_q\rest
g(N_1)\subseteq \dot{r}$ (remember that $\eta_q\in [T^q]$). Now note that 
\[p_{\eta_q\rest g(N_1)}\forces_{\p} \dot{K}\cap [N_1,g(N_1))=\emptyset.\]
Hence we easily conclude that
\[\forces_{\p}\mbox{``the increasing enumeration of $\dot{K}$ is an
unbounded real over $\V$''}\]
finishing the proof. \QED

\begin{remark}
{\em 
The forcing notion $\q_\circledcirc$ makes the ground model reals meager in a
``soft'' way: it does not add a dominating real (see \ref{forcingproperties}).
However it adds an unbounded real (just look at $\{n\in\omega: \dot{A}\cap
2^{\textstyle n}\neq\emptyset\}$, for $\dot{A}$ as in
\ref{forcingproperties}(2)). Consequently it adds a Cohen real (by
\cite{Sh:480}; note that $\q_\circledcirc$ is a Borel ccc forcing notion).
Hence we may put together \ref{unbimpcoh} and \ref{newimpunb} and we get the
following corollary.
}
\end{remark}

\begin{corollary}
Assume $\MA_\kappa(\q_\circledcirc)$. Then every ccc forcing notion of size
$\kappa$ adding new reals adds a Cohen real. \QED
\end{corollary}

\section{Anti-Martin Axiom}

In this section we are interested in axioms which are considered as strong
negations of Martin Axiom. They originated in Miller's problem if it is
consistent with $\neg\CH$ that for any ccc forcing notion of the size
$\leq\con$ there exists an $\omega_{1}$-Lusin sequence of filters (cf
\cite{MP}). The question was answered negatively by Todorcevic (cf \cite{To}).
However under some restrictions (on forcing notions and/or dense sets under
consideration) suitable axioms can be consistent with $\neg\CH$. 
These axioms were considered by van Douwen and Fleissner, who were interested
in the axiom for projective ccc forcing notions, but they needed a weakly
compact cardinal for getting the consistency (cf \cite{DF}). Cicho\'n preferred
to omit the large cardinal assumption and restricted himself to $\Sigma^1_2$
ccc forcing notions and still he was able to obtain interesting consequences
(see \cite{Ci}). Here we show how to omit the large cardinal assumption in
getting Anti--Martin Axiom for projective ccc forcing notions. This answers
Problem 6.6(2) of \cite{JR1}. 

\begin{definition}
For a forcing notion $\p$ and a cardinal $\kappa$ let $\AMA_\kappa(\p)$
be the following sentence:
\begin{quote}
{\em there exists a sequence $\langle G_i:i<\kappa\rangle$ of filters on
$\p$ such that for every maximal antichain $\A\subseteq\p$ for some
$i_0<\kappa$ we have 
\[(\forall i\geq i_0)(G_i\cap\A\neq\emptyset).\]
}
\end{quote}
For a class $\K$ of forcing notions the axiom $\AMA_\kappa(\K)$ is ``for
each $\p\in\K$, $\AMA_\kappa(\p)$ holds true''.
\end{definition}

\begin{definition}
\begin{enumerate}
\item For two models $N,M$ and an integer $n$, $M\prec_{n+1} N$ means:
\begin{quotation}
{\em
\noindent for every $\Pi_n$ formula $\varphi(x,\bar{y})$ and every
sequence $\bar{m}\subseteq M$, 

\noindent if $N\models\exists x\varphi(x,\bar{m})$ then $M\models
\exists x \varphi(x,\bar{m})$.
}
\end{quotation}
(Thus $M\prec N$ if and only if $(\forall n>0)(M\prec_n N)$.)
\item If $\p_0,\p_1$ are ccc forcing notions, $n>0$ then
$\p_0\lesdot_n\p_1$ means $\p_0\lesdot\p_1$ (i.e. $\p_0$ is a complete
suborder of $\p_1$) and
\[\forces_{\p_1}(\H(\aleph_1)^{\V[\Gamma_{\p_1}\cap\p_0]},{\in})\prec_n
(\H(\aleph_1),{\in}).\]
Instead of $\lesdot$ we may write $\lesdot_0$.
\end{enumerate}
\end{definition}

\begin{definition}
Let $\kappa$ be a cardinal number.
\begin{enumerate}
\item $\C_\kappa$ is the class of all ccc forcing notions of size
$\leq\kappa$. 
\item We inductively define subclasses $\C^n_\kappa$ of $\C_\kappa$
(for $n\leq\omega$):

$\C^0_\kappa=\C_\kappa$,

$\C_\kappa^{n+1}$ is the class of all $\p\in\C^n_\kappa$ such that for
every $\p^*\in\C^n_\kappa$
\[\p\lesdot\p^*\ \ \Rightarrow\ \ \p\lesdot_{n+1}\p^*,\]

$\C^\omega_\kappa=\bigcap\limits_{n<\omega}\C^n_\kappa$.
\end{enumerate}
\end{definition}

\begin{lemma}
\label{elecla}
Let $\kappa$ be a cardinal such that $\kappa^{\omega}=\kappa$,
$n\leq\omega$. 
\begin{enumerate}
\item If $\p_0,\p_1\in\C^n_\kappa$, $\p_0\lesdot\p_1$ then $\p_0\lesdot_n\p_1$.
\item Suppose $\delta<\kappa^+$, $\cf(\delta)\geq\omega_1$ and
$\p_i\in\C^n_\kappa$ (for $i<\delta$) are such that $i<j<\delta\ \
\Rightarrow\ \ \p_i\lesdot\p_j$. Then $\p_\delta\stackrel{\rm
def}{=}\bigcup\limits_{i<\delta}\p_i\in\C^n_\kappa$ and if $\p\in\C_\kappa$,
$\p_i\lesdot_n\p$ for every $i<\delta$ then $\p_\delta\lesdot_n\p$.
\item If $\p\in\C_\kappa$ then there is $\p^*\in\C^n_\kappa$ such that
$\p\lesdot\p^*$. 
\item If $\p\in\C_\kappa$ then there are functions
$F_k:\prod\limits_{i<\omega}\p \longrightarrow\p$ (for $k\in\omega$) such
that for every $\q\subseteq\p$:\ \ \ \ if $\q$ is closed under all $F_k$ then
$\q\lesdot_n\p$. 
\end{enumerate}
\end{lemma}

\Proof The proof is by induction on $n$. For $n=0$ there is nothing to do.
(For {\em 4.} consider functions $F_0,F_1:\prod\limits_{i\in\omega}\p
\longrightarrow\p$ such that if $\langle p_i: i\in\omega\rangle\subseteq\p$
is an antichain which is not maximal then $F_0(p_i: i<\omega)$ is a condition
incompatible with all $p_i$; if $p_i\in\p$ ($i\in\omega$) and
$p_0\compatible_{\p} p_1$ then $F_1(p_i:i\in\omega)$ is a condition stronger
than both $p_0$ and $p_1$.) So suppose that {\em 1.}--{\em 4.} hold true for
$n$ and we are proving them for $n+1$.  
\medskip

\noindent 1)\ \ \ By the definition.
\medskip

\noindent 2)\ \ \ Suppose that $\p\in\C_\kappa$, $\p_i\lesdot_{n+1}\p$ for
each $i<\delta$. By the inductive hypothesis we know that $\p_\delta\lesdot_n
\p$, $\p_\delta\in\C^n_\kappa$ and hence (by the definition of
$\C^{n+1}_\kappa$) we have $\p_i\lesdot_{n+1}\p_\delta$ for each $i<\delta$.
Suppose that $G\subseteq\p$ is a generic filter over $\V$.  Then for
$i<\delta$: 
$$(\H(\aleph_1)^{\V[G\cap\p_i]},{\in})\prec_{n+1}
(\H(\aleph_1)^{\V[G\cap\p_\delta]},{\in})\ \mbox{ and}\leqno(*)$$
$$(\H(\aleph_1)^{\V[G\cap\p_i]},{\in})\prec_{n+1}
(\H(\aleph_1)^{\V[G]},{\in}).\leqno(**)$$ 
Let $\varphi(x,\bar{y})$ be a $\Pi_n$-formula and
$\bar{y}_0\subseteq\H(\aleph_1)^{\V[G\cap\p_\delta]}$. Take $i<\delta$ 
such that $\bar{y}_0\subseteq\H(\aleph_1)^{\V[G\cap\p_i]}$ (remember
$\cf(\delta)>\omega$). If $(\H(\aleph_1^{\V[G]}),{\in})\models\exists
x\varphi(x,\bar{y}_0)$ then $(\H(\aleph_1^{\V[G\cap\p_i]}),{\in})\models\exists
x\varphi(x,\bar{y}_0)$ (by $(**)$) and
$\H(\aleph_1^{\V[G\cap\p_\delta]}),{\in})\models\exists 
x\varphi(x,\bar{y}_0)$ (by $(*)$). This shows $\p_\delta\lesdot_{n+1}\p$. To
prove that $\p_\delta\in\C^{n+1}_\kappa$ suppose that $\p\in\C^n_\kappa$,
$\p_\delta\lesdot\p$. Then for each $i<\delta$ we have $\p_i\lesdot\p$,
$\p_i\in\C^{n+1}_\kappa$ and consequently $\p_i\lesdot_{n+1}\p$. By the
previous part we get $\p_\delta\lesdot_{n+1}\p$ finishing {\em 2.}
\medskip

\noindent 3)\ \ \ Let $\p\in\C_\kappa$. By a book-keeping argument we
inductively build sequences $\langle\p_i:i\leq\kappa\rangle$ and
$\langle(p_i,\varphi_i,\dot{\tau}_i): i<\kappa\rangle$ such that for all
$i<j<\kappa$: 
\begin{enumerate}
\item $\p_i\in\C^n_\kappa$, $\p\lesdot\p_0$, $\p_i\lesdot\p_j$,
$\p_\kappa=\bigcup\limits_{i<\kappa}\p_i$,
\item $p_i\in\p_i$,
\item $\varphi_i$ is a $\Pi_n$-formula, $\dot{\tau}_i$ is a $\p_i$-name for
a finite sequence of elements of $\H(\aleph_1)$,
\item $\langle(p_i,\varphi_i,\dot{\tau}_i):\ i<\kappa\rangle$ lists all
triples $(p,\varphi,\dot{\tau})$ such that $\varphi=\varphi(x,\bar{y})$ is a
$\Pi_n$-formula, $\dot{\tau}$ is a (canonical) $\p_\kappa$-name for a finite
sequence (of a suitable length) of members of $\H(\aleph_1)$, $p\in\p_\kappa$,
\item if $i$ is limit, $\cf(i)=\omega$ then $\p_i\in\C_\kappa^n$ is such
that $\bigcup\limits_{\ell<i}\p_\ell\lesdot\p_i$,
\item if $i$ is limit, $\cf(i)>\omega$ then
$\p_i=\bigcup\limits_{\ell<i}\p_\ell\in\C^n_\kappa$,
\item {\em if\/} there is $\p^*\in\C_\kappa^n$ such that $\p_i\lesdot\p^*$
and for some $p^*\in\p^*$ we have $p^*\compatible_{\p^*}p_i$ and
\[p^*\forces_{\p^*}(\H(\aleph_1),{\in})\models\exists
x\varphi_i(x,\dot{\tau}_i)\]
{\em then\/} $\p_{i+1}$ is an example of such $\p^*$.
\end{enumerate}
The construction is fully described by the above conditions (and easy to
carry out; remember about the inductive hypothesis and the assumption that
$\kappa^\omega=\kappa$). Clearly $\p_\kappa\in\C^n_\kappa$ (by the inductive
assumption {\em 2.}) and $\p\lesdot\p_\kappa$. We have to show that actually
$\p_\kappa\in\C^{n+1}_\kappa$. Suppose not. Then we find $\p^*\in\C^n_\kappa$
such that $\p_\kappa\lesdot\p^*$ but $\p_\kappa\not\lesdot_{n+1}\p^*$. The
second means that there are a condition $p^*\in\p^*$ and a $\Pi_n$-formula
$\varphi$ and a $\p_\kappa$-name $\dot{\tau}$ for a sequence of elements of
$\H(\aleph_1)$ such that 
\[p^*\forces_{\p^*}\mbox{``}(\H(\aleph_1),{\in})\models\exists
x\varphi(x,\dot{\tau})\ \mbox{ but }\
(\H(\aleph_1)^{\V[\Gamma_{\p^*}\cap\p_\kappa]},{\in})\models\neg\exists 
x\varphi(x,\dot{\tau})\mbox{''}.\]
Take $p\in\p_\kappa$ such that $p^*\compatible_{\p^*}p$ and there is no
condition $p'\in\p_\kappa$ such that $p\leq_{\p_\kappa}p'$ and
$p'\incomp_{\p^*} p^*$. Let $i<\kappa$ be such that
$(p,\varphi,\dot{\tau})=(p_i,\varphi_i,\dot{\tau}_i)$. Condition 7 of the
construction implies that for some $p^+\in\p_{i+1}$ we have
$p^+\compatible_{\p_\kappa} p$ and
\[p^+\forces_{\p_{i+1}}(\H(\aleph_1),{\in})\models\exists
x\varphi(x,\dot{\tau}).\] 
Since $\p_{i+1}\lesdot_{n}\p_\kappa$ (the inductive hypotheses {\em 2., 1.})
we get 
\[p^+\forces_{\p_\kappa}(\H(\aleph_1),{\in})\models\exists
x\varphi(x,\dot{\tau}).\] 
The choice of $p$ implies $p^+\compatible_{\p^*}p^*$ and this provides a
contradiction as 
\[p^+\forces_{\p^*}(\H(\aleph_1)^{\V[\Gamma\cap\p_\kappa]},{\in})\models
\exists x\varphi(x,\dot{\tau}).\] 
\medskip

\noindent 4)\ \ \ Let $F^0_k:\prod\limits_{i<\omega}\p\longrightarrow\p$
(for $k\in\omega$) be functions such that if $\q\subseteq\p$, $\q$ is closed
under all $F^0_k$ then $\q\lesdot_n\p$ (they are given by the inductive
hypothesis {\em 4.}). Let $A_{i,j,k}\subseteq\omega\setminus\{0\}$ be disjoint
infinite sets (for $i,j,k\in\omega$). For a $\Pi_n$-formula
$\varphi(x,y_0,\ldots,y_{\ell-1})$ and $m\in\omega$ we choose a function
$F^{\varphi, m}:\prod\limits_{i<\omega}\p\longrightarrow\p$ satisfying the
condition described below.

\noindent Let $\langle p_m: m<\omega\rangle\subseteq\p$. For $k<\ell$ we
try to define a $\p$-name $\dot{\tau}_k$ for a real in $\baire$ by
\[(\forall m\in A_{i,j,k})(p_m\forces_{\p}\dot{\tau}_k(i)=j).\]
If this definition is correct then we ask if these reals encode (in the
canonical way) elements of $\H(\aleph_1)$ (which we identify with the names
$\dot{\tau}_k$ themselves). If yes then we ask if 
\[p_0\forces_{\p}(\H(\aleph_1),{\in})\models\exists
x\varphi(x,\dot{\tau}_0,\ldots,\dot{\tau}_{\ell-1}).\]
If the answer is positive then we fix a $\p$-name $\dot{\tau}$ for (a real
encoding) a member of $\H(\aleph_1)$ such that
\[p_0\forces_{\p}(\H(\aleph_1),{\in})\models
\varphi(\dot{\tau},\dot{\tau}_0,\ldots,\dot{\tau}_{\ell-1}).\]
This name can be represented similarly as names $\dot{\tau}_k$ (for
$k<\ell$) so we have a sequence $\langle q_m: m<\omega\rangle\subseteq\p$
encoding it. Finally we want $F^{\varphi,m}$ to be such that if the above
procedure for $\langle p_m:m<\omega\rangle$ works then
$F^{\varphi,m}(p_m:m<\omega)=q_m$. 

\noindent Now take all the functions $F^0_k$, $F^{\varphi,m}$; it is easy to
check that they work.
\medskip

Lastly note that the case $n=\omega$ follows immediately from the lemma
for $n<\omega$. (For {\em 3.} construct an increasing sequence
$\langle\p_i: i<\omega_1\rangle$ such that $\p\lesdot\p_0$ and if
$\lambda<\omega_1$ is limit, $k<\omega$ then
$\p_{\lambda+k}\in\C^k_\kappa$.) \QED

\begin{theorem}
\label{getama}
\hspace{0.15in}
Suppose that $\theta,\kappa$ are cardinals such that
$\aleph_1\leq\theta=\cf(\theta)\leq\kappa=\kappa^{\omega}$. Then there
exists a ccc forcing notion $\p$ such that 
\[\forces_{\p}\con=\kappa\ \&\ \AMA_\theta(\mbox{projective ccc}).\]
\end{theorem}

\Proof The forcing notion $\p$ which we are going to construct will be
essentially a finite support iteration of length $\kappa\cdot\theta$ of ccc
forcing notions. One could try to force with ``all possible ccc orders'' in the
iteration. However some care is necessary to make sure that several notions
(including ``being a maximal antichain'') are sufficiently absolute for
intermediate stages. Therefore we use forcing notions from the class
$\C^\omega_\kappa$. So we inductively build sequences $\langle\p_i:
i\leq\kappa\cdot\theta\rangle$ and $\langle(\varphi_i,\psi_i,\dot{\tau}_i):
i<\kappa\cdot\theta\rangle$ such that for all $i<j<\kappa\cdot\theta$:
\begin{enumerate}
\item $\p_i\in\C^\omega_\kappa$, $\p_{\kappa\cdot\theta}=
\bigcup\limits_{i<\kappa\cdot\theta}\p_i\in\C^\omega_\kappa$,
\item $\p_i\lesdot\p_j$,
\item $\langle(\varphi_i,\psi_i,\dot{\tau}_i): i<\kappa\cdot\theta\rangle$
lists with cofinal repetitions all triples $(\varphi,\psi,\dot{\tau})$ such
that $\varphi$ is a formula with $n+1$ variables, $\psi$ is a formula with
$n+2$ variables and $\dot{\tau}$ is a $\p_{\kappa\cdot\theta}$-name for a
sequence of length $n$ of elements of $\H(\aleph_1)$,
\item if $\dot{\tau}_i$ is a $\p_i$-name and
\[\forces_{\p_i}\mbox{``}\langle\varphi_i(x,\dot{\tau}_i),
\psi_i(x_0,x_1,\dot{\tau}_i)\rangle \mbox{ defines in } (\H(\aleph_1),{\in})
\mbox{ a ccc partial order }\dot{\q}_i\mbox{''}\]
then $\p_i*\dot{\q}_i\lesdot\p_{i+1}$.
\end{enumerate}
It is easy to carry the construction (use a book-keeping argument, remembering
$\kappa^\omega=\kappa$ plus lemma~\ref{elecla}). We want to show that
$\p=\p_{\kappa\cdot\theta}$ has the required properties. Easily
$\forces_{\p}\con=\kappa$. Now suppose that $G\subseteq\p$ is a generic filter
over $\V$ and work in $\V[G]$. 
\medskip

Assume that $\q$ is a projective ccc forcing notion and thus it is definable
in $(\H(\aleph_1),{\in})$. Thus we have formulas $\varphi(x,\bar{y})$ and
$\psi(x_0,x_1,\bar{y})$ and a sequence $\bar{r}\subseteq\H(\aleph_1)$ such
that
\[\q=\{x\in\H(\aleph_1): (\H(\aleph_1),{\in})\models\varphi(x,\bar{r})\}\]
\[\leq_{\q}=\{(x_0,x_1)\in\H(\aleph_1)\times\H(\aleph_1):
(\H(\aleph_1),{\in})\models\psi(x_0,x_1,\bar{r})\}.\] 
Let $\dot{\tau}$ be a $\p$-name for $\bar{r}$. We may assume that 
\[\forces_{\p}\mbox{``}\langle\varphi(x,\dot{\tau}),\psi(x_0,x_1,\dot{\tau})
\rangle \mbox{ defines (in $(\H(\aleph_1),{\in})$) a ccc partial order''}.\]
There is an increasing cofinal in $\kappa\cdot\theta$ sequence $\langle i_j:
j<\theta\rangle$ such that $\dot{\tau}$ is a $\p_{i_0}$-name and
$(\varphi_{i_j},\psi_{i_j},\dot{\tau}_{i_j})= (\varphi,\psi,\dot{\tau})$.
Since $\p,\p_{i_j}\in\C^\omega_\kappa$ we have that 
\[\forces_{\p}\mbox{``}(\H(\aleph_1)^{\V[\Gamma_{\p}\cap\p_{i_j}]},{\in})\prec
(\H(\aleph_1),{\in})\mbox{''}\] 
and hence the formulas $\langle\varphi(x,r),\psi(x_0,x_1,r)\rangle$ define (in
$(\H(\aleph_1)^{\V[G\cap\p_{i_j}]},{\in})$) the partial order
$\q\cap\H(\aleph_1)^{\V[G\cap\p_{i_j}]}$. The incompatibility relation in this
partial order is expressible in $(\H(\aleph_1),{\in})$ and thus it is the
restriction of $\incomp_{\q}$. Consequently $\q\cap\H(\aleph_1)^{\V[G\cap
\p_{i_j}]}$ is ccc in $\V[G]$ and hence in $\V[G\cap\p_{i_j}]$. Hence in
$\V[G\cap\p_{i_{j+1}}]$ we have a filter $G^*_j\subseteq\q\cap
\H(\aleph_1)^{\V[G\cap\p_{i_j}]}$ generic over $\V[G\cap\p_{i_j}]$ (here we
apply condition 4 of the construction). Look at the sequence $\langle G^*_j:
j<\theta\rangle$. Let $\A\subseteq\q$ be a maximal antichain. It is countable
and hence for sufficiently large $j<\theta$ we have $\A\in\V[G\cap\p_{i_j}]$.
Moreover the antichain can be coded as a one real and the fact that it is a
maximal antichain in the partial order defined by $\langle\varphi,\psi\rangle$
is expressible in $(\H(\aleph_1),{\in})$. Applying $\p_{i_j}\in
\C^\omega_\kappa$ we get that  
\[\V[G\cap\p_{i_j}]\models\A\mbox{ is a maximal antichain in }\q\cap
\H(\aleph_1)^{\V[G\cap\p_{i_j}]}.\] 
Consequently for sufficiently large $j<\theta$ we have
\[G^*_j\cap\A\neq\emptyset.\]
This finishes the proof.\QED

\begin{remark}
{\em 
In \ref{elecla} and \ref{getama} we used $\H(\aleph_1)$ as we were mainly
interested in $\AMA_{\omega_1}$ and projective ccc forcing notions. But we may
replace it by $\H(\chi)$ for any uncountable regular cardinal $\chi$ such that
$\sum\limits_{\alpha<\chi}\kappa^{|\alpha|}=\kappa$. Then in \ref{elecla}(2)
we consider $\delta<\kappa^+$ such that $\cf(\delta)\geq \chi$ and in
\ref{getama} we additionally assume that $\theta\geq\chi$. 
}
\end{remark}

\section{Absoluteness and embeddings}
In this section we answer positively Problem 4.4 of \cite{JR1} (see
\ref{getabs}) and we give a negative answer to Problem 3.3 of \cite{JR1} (see
\ref{ad3.3con}). 

\begin{definition}
Let $\p$ be a forcing notion. We say that $\Sigma^1_n(\p)$-absoluteness
holds if for every $\Sigma^1_n$ formula $\varphi$ (with parameters in $\V$)
and a generic filter $G\subseteq\p$ over $\V$
\begin{quote}
\noindent $\V[G]\models\varphi$ if and only if $\V\models\varphi$.
\end{quote}
\end{definition}
Obviously $\Sigma^1_2(\p)$-absoluteness holds for any forcing notion $\p$.

\begin{theorem}
\label{getabs}
Assume $\MA_{\omega_1}(\p)$. Then $\Sigma^1_3(\p)$--absoluteness holds.
\end{theorem}

\Proof Suppose that $\varphi$ is a $\Sigma^1_3$ sentence (with a parameter
$a\in\baire$). Using the tree representation of $\Pi^1_2$--sets we find
a tree $T$ (constructible from $a$) over $\omega\times\omega_1$ such
that
\[\begin{array}{ll}
\varphi\equiv & (\exists x\in\baire)(\forall f\in{\omega_1}^{\textstyle
\omega})(\exists n\in\omega)(\langle x\restriction n, f\restriction
n\rangle\notin T)\\
\equiv & (\exists x\in\baire)(\mbox{the tree $T(x)$ is well founded}).
\end{array}\]
(For $x\in\baire$, $T(x)$ is the tree on $\omega_1$ consisting of all
$\bar{\alpha}\in{\omega_1}^{\textstyle {<}\omega}$ such that $\langle
x\restriction \lh(\bar{\alpha}),\bar{\alpha}\rangle\in T$.) Moreover, as by
$\MA_{\omega_1}(\p)$ we know that $\forces_{\p}\omega_1^{\V}=\omega_1$, the
tree $T$ represents $\varphi$ in $\V^{\p}$ too:
\[\forces_{\p}\mbox{``}\varphi\equiv (\exists x\in\baire)(\mbox{the tree
$T(x)$ is well founded})\mbox{''}. \]
Suppose now that $\forces_{\p}\varphi$. Then we have a $\p$-name $\dot{r}$ for
a real in $\baire$ such that 
\[\forces_{\p}\mbox{``the tree $T(\dot{r})$ is well founded''}.\]
Consequently we have a $\p$-name $\dot{\rho}$ for a function such that 
\[\forces_{\p}\mbox{``}\dot{\rho}:T(\dot{r})\longrightarrow{\rm Ord}\ \mbox{
is a rank function''.}\]
For $n\in\omega$, $\bar{\alpha}\in{\omega_1}^{\textstyle n}$ put
\[J^0_n=\{p\in\p: (\exists m\in\omega)(p\forces_{\p}\dot{r}(n)=m)\},\]
\[\begin{array}{ll}
J^1_{\bar{\alpha}}=\{p\in\p: &\mbox{either }
p\forces_{\p}\bar{\alpha}\notin T(\dot{r})\\
\  &\mbox{or } (\exists\xi\in {\rm Ord})(p\forces_{\p}\bar{\alpha}\in
T(\dot{r})\ \&\ \dot{\rho}(\dot{r}\restriction n,\bar{\alpha})=\xi)\}.\\
\end{array}\]
Clearly these are dense subsets of $\p$. By $\MA_{\omega_1}(\p)$ we
find a filter $G$ on $\p$ such that $G\cap J^0_n\neq\emptyset$ for
$n\in\omega$ and $G\cap J^1_{\bar{\alpha}}\neq\emptyset$ for
$\alpha\in{\omega_1}^{\textstyle{<} \omega}$. Using this filter we may
interpret the name $\dot{r}$ to get $r=\dot{r}^G\in\baire$. Moreover we
may interpret the name $\dot{\rho}$ to get a function $\rho=\rho^G:
T(r)\longrightarrow {\rm Ord}$: 
\[\rho(r\restriction n,\bar{\alpha})=\xi\ \mbox{ iff }\ (\exists p\in
G)(p\forces_{\p}\bar{\alpha}\in T(\dot{r})\ \&\
\dot{\rho}(\dot{r}\restriction n,\bar{\alpha})=\xi).\]
[Note that this really defines a function from $T(r)$ to ordinals: suppose
that $\langle r\restriction n, \bar{\alpha}\rangle\in T$. First we find $p\in
G\cap\bigcap\limits_{m<n}J^0_m$; then clearly 
\[p\forces_{\p}\mbox{``}\dot{r}\restriction n = r\restriction n\ \mbox{
and }\ \bar{\alpha}\in T(\dot{r})\mbox{''.}\]
Thus if $p'\in J^1_{\bar{\alpha}}\cap G$ then $p'\not\forces_{\p}\bar{\alpha}
\notin T(\dot{r})$ and hence for some ordinal $\xi$ we have $p'\forces_{\p}
\bar{\alpha}\in T(\dot{r})\ \&\ \dot{\rho}(\dot{r}\restriction n,\bar{\alpha})
=\xi$. Moreover if $\xi_0,\xi_1$ are such that for some $p^0,p^1\in G$ we have 
\[p^i\forces_{\p}\bar{\alpha}\in T(\dot{r})\ \&\
\dot{\rho}(\dot{r}\restriction n,\bar{\alpha})=\xi_i\]
then (as $p^0\compatible_{\p} p^1$) we cannot have $\xi_0\neq\xi_1$.]

\noindent We claim that $\rho$ is a rank function on $T(r)$. Suppose that
$n_0<n_1$, $\bar{\alpha}_0\in{\omega_1}^{\textstyle n_0}$, $\bar{\alpha}_1\in
{\omega_1}^{\textstyle n_1}$, $\bar{\alpha}_0\subsetneq \bar{\alpha}_1$ and
$\langle r\restriction n_0,\bar{\alpha}_0\rangle, \langle r\restriction
n_1,\bar{\alpha}_1\rangle\in T$. Take a condition $p\in G\cap\bigcap
\limits_{m<n_1} J^0_m$. Then  
\[p\forces_{\p}\mbox{``}\dot{r}\restriction n_1=r\restriction n_1\ \&\
\bar{\alpha}_0,\bar{\alpha}_1\in T(\dot{r})\mbox{''}.\]
Next choose conditions $p^0,p^1\in G$ such that
\[p^i\forces_{\p}\mbox{``}\bar{\alpha}_i\in T(\dot{r})\ \&\
\dot{\rho}(\dot{r}\restriction n_i,\bar{\alpha}_i)=\rho(\bar{\alpha}_i)
\mbox{''}.\]
Take $p^*\in G$ stronger than $p^0,p^1,p$. Since $\dot{\rho}$ is (forced to
be) a rank function on $T(\dot{r})$ we have 
\[p^*\forces_{\p}\rho(\bar{\alpha}_0)=\dot{\rho}(\dot{r}\restriction
n_0,\bar{\alpha}_0)> \dot{\rho}(\dot{r}\restriction
n_1,\bar{\alpha}_1)=\rho(\bar{\alpha}_1).\] 
Hence $\rho(\bar{\alpha}_0)>\rho(\bar{\alpha}_1)$ and we may conclude our
theorem: the tree $T(r)$ is well founded so $\V\models\varphi$.\QED

\begin{proposition}
Suppose that $\q$ is a ccc Souslin forcing notion (i.e. $\q$, $\leq_{\q}$ and
$\incomp_{\q}$ are $\Sigma^1_1$--sets), $\dot{r}$ is a $\q$--name for a
function from $\fs$ to $2$. Let 
\[A[\dot{r}]\stackrel{\rm def}{=}\{\eta\in\can: (\exists
p\in\q)(\forall^\infty
m\in\omega)(p\forces_{\q}\mbox{``}\dot{r}(\eta\restriction
m)=1\mbox{''})\}.\]
Then $A[\dot{r}]$ is an analytic set.
\end{proposition}

\Proof For each $\nu\in\fs$ choose a maximal antichain $\langle p_{\nu,l}:
l\in\omega\rangle$ in $\q$ and a set $I_\nu\subseteq\omega$ such that for
each $l\in\omega$:
\[l\in I_\nu\ \Rightarrow\ p_{\nu,l}\forces_{\q} \dot{r}(\nu)=0\ \ \mbox{
and }\ \ l\notin I_\nu\ \Rightarrow\ p_{\nu,l}\forces_{\q} \dot{r}(\nu)=1.\]
Now note that for each $\eta\in\can$ we have
\[\eta\in A[\dot{r}]\equiv (\exists p\in\q)(\forall^\infty n\in\omega)(\forall
l\in I_{\eta\restriction n})(p_{\eta\restriction n,l}\incomp_{\q} p).\] 
 \QED

\begin{proposition}
For every $A\subseteq\can$ there exist a $\sigma$-centered forcing notion
$\q^A$ and a $\q^A$--name $\dot{r}$ (for a function from $\fs$ to $2$) such
that $A=A[\dot{r}]$ and $\|\q\|=\|A\|+\aleph_0$. 
\end{proposition}

\Proof The forcing notion $\q^A$ is defined by
\medskip

\noindent{\bf conditions} are pairs $\langle r,w\rangle$ such that $r$ is a
finite function, $\dom(r)\subseteq\fs$, $\rng(r)\subseteq 2$ and $w\in
[A]^{\textstyle <\!\omega}$,

\noindent{\bf the order} is such that $\langle r_1,w_1\rangle\leq \langle
r_2,w_2\rangle$ if and only if $r_1\subseteq r_2$, $w_1\subseteq w_2$ and
\[(\forall\nu\in\dom(r_2)\setminus\dom(r_1))([(\exists\eta\in w_1)(\nu
\subseteq\eta)]\ \Rightarrow\ r_2(\nu)=1).\]
\medskip

\noindent The $\q^A$--name $\dot{r}$ is such that 
\[\forces_{\q^A}\dot{r}=\bigcup\{r: (\exists w)(\langle
r,w\rangle\in\Gamma_{\q^A})\}.\] 
It should be clear that $\q^A$ is $\sigma$-centered,
$\|\q^A\|=\|A\|+\aleph_0$ and 
\[\forces_{\q^A}\dot{r}:\fs\longrightarrow 2.\]
Moreover for each $\eta\in\can$ and $\langle r,w,\rangle\in \q^A$:
\[(\forall^\infty m)(\langle
r,w\rangle\forces_{\q^A}\dot{r}(\eta\restriction m)=1)\ \ \mbox{ iff }\ \
\eta\in w.\]
Consequently $A=A[\dot{r}]$.\QED

\begin{corollary}
\label{ad3.3con}
If $A\subseteq\can$ is not analytic then $\q^A$ cannot be completely
embedded into a ccc Souslin forcing. In particular, if $\con>\aleph_1$ then
there is a $\sigma$-centered forcing notion of size $\aleph_1$ which cannot
be completely embedded into a ccc Souslin forcing notion.\QED
\end{corollary}


\begin{thebibliography}{XXSh000}

\bibitem[BaJu]{BaJu} Bartoszy\'nski, Tomek and Judah, Haim, {\bf Set Theory: On
the Structure of the Real Line}, A K Peters, Wellesley, Massachusetts, 1995.

\bibitem[Ci]{Ci} Cicho\'n, Jacek, {\em Anti--Martin Axiom}, circulated notes
(1989). 

\bibitem[DF]{DF} van~Douwen, Eric K. and Fleissner, William G., {\em Definable
Forcing Axiom: An Alternative to Martin's Axiom}, Topology and its
Applications vol.35(1990): 277--289.

\bibitem[Je]{Je} Jech, Thomas, {\bf Set Theory}, Academic Press 1978.

\bibitem[JR1]{JR1} Judah, Haim and Ros{\l}anowski, Andrzej, {\em Martin
Axiom and the continuum}, {\bf Journal of Symbolic Logic}, vol.60(1995):
374--391. 

\bibitem[MP]{MP} Miller, Arnold and Prikry, Karel, {\em When the continuum
has cofinality $\omega_{1}$}, Pacific Journal of Mathematics,
vol.115(1984): 399-407.

\bibitem[Sh:480]{Sh:480} Shelah, Saharon, {\em How special are Cohen and
Random forcings}, Israel Journal of Mathematics, vol.88(1994): 153--174.

\bibitem[Sh:f]{Sh:f} Shelah, Saharon, {\bf Proper and improper forcing},
Perspectives in Mathematical Logic, Springer, accepted.

\bibitem[To]{To} Todorcevic, Stevo, {\em Remarks on Martin's Axiom and the
Continuum Hypothesis}, Canadian Journal of Mathematics, vol.43 (1991):
832--851. 

\end{thebibliography}
\end{document}